\documentclass[11pt]{amsart}
\usepackage{amssymb, latexsym}
\theoremstyle{plain}
\newtheorem{theorem}{Theorem}
\newtheorem{corollary}{Corollary}

\newtheorem{proposition}{Proposition}
\newtheorem{assumption}{Assumption}
\newtheorem*{EP1}{Theorem EP1}
\newtheorem*{EP2}{Theorem EP2}
\newtheorem*{EP3}{Theorem EP3}
\newtheorem*{EP4}{Theorem EP4}
\newtheorem*{p}{Theorem BP}
\newtheorem*{compare}{Comparison Result}

\theoremstyle{remark}

\newtheorem*{remarK 1}{Remark 1}
\newtheorem*{remarK 2}{Remark 2}
\newtheorem*{Remark 1}{Remark 1}
\newtheorem*{Remark 2}{Remark 2}
\newtheorem*{Remark 3}{Remark 3}
\newtheorem*{remark 1}{Remark 1}
\newtheorem*{remark 2}{Remark 2}
\newtheorem*{remark 3}{Remark 3}
\newtheorem*{remark 4}{Remark 4}

\newtheorem*{REmark 1}{Remark 1}
\newtheorem*{REmark 2}{Remark 2}
\newtheorem*{REmark 3}{Remark 3}

\newtheorem*{remark}{Remark}

\theoremstyle{definition}

\newtheorem*{definition}{Definition}

\numberwithin{equation}{section}

\begin{document}

\title {The compact support property for measure-valued processes}

\author{J{\'a}nos Engl\"ander\ and\ Ross G. Pinsky}
\address{Department of Statistics and Applied Probability\\
University of California, Santa Barbara, CA 93106-3110, USA
 \\
and Department of Mathematics\\
Technion---Israel Institute of Technology\\ Haifa, 32000\\ Israel}

\email{englander@pstat.ucsb.edu; pinsky@math.technion.ac.il}

\urladdr{http://www.pstat.ucsb.edu/faculty/englander;
www.math.technion.ac.il/people/pinsky/}



\date{}

\begin{abstract}

The purpose of this article is to give a rather thorough
understanding
of the compact support property for measure-valued
processes corresponding to semi-linear equations of the
form
\[
\begin{aligned}& u_t=Lu+\beta u-\alpha u^p \ \ \text{in} \   R^d\times (0,\infty),
\ p\in(1,2];\\
&u(x,0)=f(x) \ \ \text{in}\ R^d;\\
&u(x,t)\ge0 \ \ \text{in} \ R^d\times[0,\infty).
\end{aligned}
\]
In particular, we shall investigate how the interplay
between the underlying motion (the diffusion process corresponding
to $L$)  and the branching  affects
the compact support property.
In \cite{EP99}, the compact support property
was shown to be equivalent to a certain analytic criterion concerning
uniqueness of the Cauchy problem for the semi-linear parabolic
equation related to the measured valued process.
In  a subsequent paper \cite{EP03}, this analytic property was investigated
purely from the point of view of partial differential equations.
Some of the results obtained in this latter paper
yield  interesting results concerning the compact support property.
In this paper, the results from \cite{EP03} that are relevant
to the compact support property are presented, sometimes  with extensions.
These results are interwoven with   new results and
some informal heuristics.
Taken together, they yield a  rather comprehensive picture of the compact support
property.
\it Inter alia\rm, we show that the concept of a measure-valued
process \it hitting\rm\ a point can be investigated via
the compact support property, and suggest an alternate
proof of a result concerning the hitting of points by super-Brownian
motion.

\end{abstract}

\maketitle

\section{Introduction and Statement of Results}\label{S:intro}

The purpose of this article is to give a rather thorough
understanding
of the compact support property for measure-valued
diffusion processes.
In particular, we shall investigate how the interplay
between the underlying motion and the branching affects
the compact support property.
In \cite{EP99}, the compact support property
was shown to be equivalent to a certain analytic criterion concerning
uniqueness of the Cauchy problem for the semi-linear parabolic
equation related to the measured valued process.
In  a subsequent paper \cite{EP03}, this analytic property was investigated
purely from the point of view of partial differential equations.
Some of the results obtained in this latter paper
yield  interesting results concerning the compact support property.
In this paper, the results from \cite{EP03} that are relevant
to the compact support property are presented, sometimes  with extensions.
These results are interwoven with   new results and
some informal heuristics.
Taken together, they yield a  rather comprehensive picture of the compact support
property.
\it Inter alia\rm, we show that the concept of a measure-valued
process \it hitting\rm\ a point can be investigated via
the compact support property
and suggest an alternate
proof of a result concerning the hitting of points by super-Brownian
motion.

We state all of our results for the case that the underlying space is $R^d$; however, all the
results also hold for generic domains in $R^d$---see \cite{EP99}, where
 the superprocesses studied in this paper are constructed on generic Euclidean domains.

We begin by defining the measure-valued processes under study.
Let
$$
L=\frac12\sum_{i,j=1}^d a_{i,j}\frac{\partial^2}{\partial x_i\partial
x_j}+\sum_{i=1}^d b_i\frac\partial{\partial x_i},
$$
where $a_{i,j},b_i\in C^{\alpha}(R^d)$, with
$\alpha\in(0,1)$, and $\{a_{i,j}\}$ is
strictly elliptic; that is,
$\sum_{i,j=1}^da_{i,j}(x)\nu_i\nu_j>0$, for all $x\in R^d$ and
$\nu\in R^d-\{0\}$. Let $Y(t)$ denote the diffusion process
corresponding to the generalized martingale problem for $L$ on
$R^d$ \cite{P95}, and denote corresponding probabilities by
$\mathcal P_\cdot$. Denote the lifetime of $Y(t)$ by
$\tau_\infty$. One has $\tau_\infty=\lim_{n\to\infty}\tau_n$,
where $\tau_n=\inf\{t\ge0:|Y(t)|\ge n\}$. Recall that $Y(t)$ is
called \it non-explosive\rm\ (or conservative) if $\mathcal
P_x(\tau_\infty<\infty)=0$, for some, or equivalently all, $x\in
R^d$; otherwise $Y(t)$ is called \it explosive.\rm\ The process
$Y(t)$  serves as the underlying motion of the measure-valued
process.

The branching mechanism is of the form
$\Phi(x,z)= \beta(x)z-\alpha(x)z^p$, where $p\in(1,2]$,
$\beta$ is  bounded from above, $\alpha>0$,
and $\alpha, \beta\in C^\kappa(R^d)$, for some $\kappa\in (0,1]$.
A finite measure-valued
process $X(t)=X(t,\cdot)$ is   the Markov process defined uniquely
via  the following log-Laplace equation:
\begin{equation}\label{Def:MVD}
E(\exp(-<f, X(t)>)|X(0)=\mu)=\exp(-<u_f(\cdot,t),\mu>),
\end{equation}
for $f\in C^+_c(R^d)$,
the space of compactly supported, nonnegative, continuous
functions on $R^d$, and for finite initial measures
$\mu$, where
$u_f$ is the minimal positive solution to the evolution equation

\begin{equation}\label{evoequ}
\begin{aligned}& u_t=Lu+\beta u-\alpha u^p \ \ \text{in} \   R^d\times (0,\infty);\\
&u(x,0)=f(x) \ \ \text{in}\ R^d;\\
&u(x,t)\ge0 \ \ \text{in} \ R^d\times[0,\infty).
\end{aligned}
\end{equation}
(For the construction, see \cite{EP99}.)
The measure for the process started from $\mu$ will be denoted
by $P_\mu$, and its expectation operator will be denoted by $E_\mu$.

We recall the compact support property.

\begin{definition}   Let $\mu\in \mathcal M_F(R^d)$ be compactly
supported. The measure-valued process corresponding
to $P_\mu$  possesses the  compact support property
if
\begin{equation}
P_\mu(\bigcup_{0\le s\le t} \text{supp}\ \  X(s)\ \text{ is bounded})=1, \
\text{for all}\ t\ge0.
\end{equation}
\end{definition}

\begin{remark} The parameter $\beta$ may be thought of as the mass
creation parameter (see  the discussion of the particle process
approximation to the measure-valued process  at
the end of this section).
It is possible to extend
the construction of the measure-valued process    to
certain $\beta$ which are unbounded from above; namely, to
those $\beta$  for which the generalized principal eigenvalue of the
operator $L+\beta$ is finite. However, the resulting process
has paths which are mutually absolutely continuous with respect to
another measure-valued process whose mass creation parameter
$\beta$ \it is\rm\ bounded from above \cite{EP99}. Thus, the compact support
property will hold for the former process
if and only if it holds for the
latter one.
 \it Without further mention,
it will  always be assumed in this paper
that $\beta$ is bounded from above.\rm\
\end{remark}

There are four objects,
corresponding to four different underlying probabilistic effects,
 which can influence  the compact support
property:
\begin{enumerate}
\item
$L$, the operator  corresponding to the underlying motion;
\item
$\beta$, the mass creation parameter of the branching mechansim;
\item $\alpha$, the nonlinear component
of the branching mechanism, which can be thought of as the variance
parameter if $p=2$;
\item
$p$, the power of the nonlinearity, which
is  the scaling power and is connected to the fractional  moments
of the offspring distribution in the particle process
approximation to the measure-valued process.

\noindent  (For (2), (3) and (4) above, see
the discussion of the particle process approximation to the measure-valued
process at the end of this section.)

\end{enumerate}

We shall see that both $L$ and $\alpha$ play a large role in
determining whether or not the compact support property holds;
$\beta$ and $p$ play only a
 minor role.

In \cite{EP99}, the compact support property was shown to be equivalent
to a uniqueness property for solutions to \eqref{evoequ}.
\begin{EP1}
The compact support property holds
for one, or equivalently all, nonzero, compactly supported initial measures
$\mu$
if and only if there
are no nontrivial solutions
to \eqref{evoequ} with initial data $f\equiv0$.
\end{EP1}

\begin{Remark 1}
We emphasize that
the uniqueness property in
 Theorem EP1 concerns \it all\rm\ classical solutions to
\eqref{evoequ},
with no growth restrictions. If one restricts to mild solutions---solutions
which solve an integral equation involving the linear semigroup corresponding
to the operator $L$---then, for example, uniqueness holds in this class
if $\alpha, \beta$ and the coefficients of $L$ are bounded \cite{P83};
yet, these conditions certainly
do not guarantee uniqueness in the
class of all positive, classical solutions (see Theorem \ref{smallalpha}
below).

\end{Remark 1}
\begin{Remark 2}
In fact, the proof of Theorem EP1 shows that
there exists a maximal solution $u_{max}$ to \eqref{evoequ} with
initial data $f=0$, and
\[
P_\mu(\bigcup_{0\le s\le t}\text{supp}\ \  X(s)\ \text{ is bounded})=
\exp(-\int_{R^d}u_{max}(x,t)\mu(dx)),
\]
for compactly supported $\mu$.
This shows that
when the compact support property fails,
 the onset of the failure
is gradual; that is, as a function of $t$,\linebreak
$P_\mu(\bigcup_{0\le s\le t}\text{supp}\  X(s)\ \text{ is bounded})$
is continuous and equal to 1 at time $t=0$.
This behavior is in contrast to the behavior
of the measure-valued process corresponding to the semi-linear operator
$u_t=\Delta u-u\log u$,
investigated recently in \cite{FS04}. This process is obtained
as a weak limit as $p\to1$  of the processes corresponding to
the semi-linear operators $u_t=\Delta u+\frac1{p-1}u-\frac1{p-1}u^p$.
Unlike the measure-valued processes defined above,
this process is  immortal; that is, $P_\mu(X(t)=0)=0$, for all
$t\ge0$. Furthermore, it is shown that
$P_\mu(\bigcup_{0\le s\le t}\text{supp}\ \  X(s)\ \text{ is bounded})=0$,
for all $t>0$.
Thus,  the onset of the  failure of the compact support property
is instantaneous.
Theorem EP1 is not valid for this process. Indeed, the proof of Theorem EP1
requires  the fact that a maximal solution exists for \eqref{evoequ}
with initial condition $f=0$.
The existence of such a maximal solution
is essentially  equivalent to the existence of a universal, a priori upper bound
on all solutions to \eqref{evoequ}; that is, the existence of  a finite function
$M(x,t)$ on $R^d\times(0,\infty)$ such that $u(x,t)\le M(x,t)$,
for all $(x,t)\in R^d\times(0,\infty)$, and all solutions $u$ to \eqref{evoequ}.
Such a universal a priori upper bound does not exist for the equation
$u_t=\Delta u-u\log u$.
In \cite{P05},   a more or less necessary and sufficient
condition on the nonlinear term
(independent of the operator $L$)
is given for the existence of such a bound.

\end{Remark 2}

\begin{Remark 3} Theorem EP1 suggests a parallel between the compact support
property for measure-valued processes and the non-explosion
property for diffusion processes.
Indeed, the non-explosion property for the  diffusion process
$Y(t)$ corresponding to the operator $L$
is equivalent to the nonexistence of nontrivial, \it bounded\rm\, positive solutions
to the linear Cauchy problem with 0 initial data:
\begin{equation}\label{linequ}
\begin{aligned}&u_t=Lu\ \ \text{in}\ R^d\times(0,\infty);\\
&u(x,0)=0\ \ \text{in}\ R^d;\\
&u\ge0 \ \ \text{in} \ R^d\times(0,\infty).
\end{aligned}
\end{equation}
(For one direction of this result, note that
$u(x,t)\equiv\mathcal P_x(\tau_\infty\le t)$ serves as a nontrivial solution to \eqref{linequ}
in the explosive case.)
It is natural for bounded, positive solutions to be the relevant  class of solutions
in the linear case and for positive solutions to be the relevant class
of solutions in the semi-linear case. Indeed, by Ito's formula,
the  probabilities for certain events related to $Y(t)$ are obtained as
bounded, positive solutions to the linear equation,
and by the log-Laplace equation, the negative of the logarithm of the probability
of certain events related to $X(t)$ can be obtained as positive solutions
to the semi-linear equation.
\end{Remark 3}
The class of operators $L$ satisfying the following assumption
will play an important role.
\begin{assumption}\label{quadlin} For some $C>0$,\hfill
\begin{enumerate}
\item $\sum_{i,j=1}^na_{ij}(x)\nu_i\nu_j
\le C|\nu|^2(1+|x|^2), \ x,\nu\in R^d$;
\item $|b(x)|\le C(1+|x|),\  x\in R^d$.
\end{enumerate}
\end{assumption}

The next theorem culls some results from \cite{EP03}
and applies them to
the probabilistic setting at hand.

\begin{EP2} Let $p\in(1,2]$ and let the coefficients of $L$ satisfy
Assumption \ref{quadlin}. \hfill
\begin{enumerate}
\item There is no nontrivial solution to \eqref{linequ};
thus, the diffusion process $Y(t)$ does not explode.
\item If
\[\inf_{x\in R^d}\alpha(x)>0,
\]
then
there is no nontrivial solution to \eqref{evoequ} with initial data
$f=0$; thus,
the compact support property holds for $X(t)$.
\end{enumerate}
\end{EP2}
\begin{proof}
[Proof of part (1)]
For the proof that there are no nontrivial solutions to \eqref{linequ}
if Assumption \ref{quadlin} holds, see
 \cite[Proposition 5 and Remark 1 following it]{EP03}.
Non-explosiveness of $Y(t)$ then follows from the parenthetical
sentence following \eqref{linequ} above.

\it \noindent \it Proof of  part (2).\rm\ For the proof that there are no nontrivial solutions
to \eqref{evoequ}, see    \cite[Theorem 2]{EP03}. The fact that the compact
support property holds then follows from Theorem \nolinebreak EP1 above.
\end{proof}

The conditions in Assumption 1 are classical conditions which
arise frequently in the theory of diffusion processes.
Theorem EP2 shows that if the coefficients  of $L$ obey this condition
and if the branching coefficient $\alpha$ is bounded away from zero, then
everything is well behaved---neither can the underlying
diffusion process explode nor can the measure-valued process
fail to possess the compact support property.

The following result
shows that the compact support property can fail if
$\inf_{x\in R^d}\alpha(x)=0$. It also
demonstrates that  the
effect of $\alpha$ on the compact support property
cannot be studied in isolation, but in fact depends on the
underlying diffusion.

\begin{theorem}\label{smallalpha}
Let $p\in(1,2]$ and let
\[
\begin{aligned}
&L=\frac12\sum_{i,j=1}^da_{i,j}\frac{\partial^2}{\partial x_i\partial x_j},\
\text{where} \\
& C_0^{-1}(1+|x|)^m\le \sum_{i,j=1}^da_{i,j}(x)
\le C_0(1+|x|)^m,\ m\in[0,2],\
 \text{for some}\ C_0>0.
\end{aligned}
\]
\begin{enumerate}
\item If
\[
\alpha(x)\ge C_1\exp(-C_2|x|^{2-m}),
\]
for some $C_1,C_2>0$,
then  the compact support property holds for $X(t)$.
\item If
\[
\alpha(x)\le C\exp(-|x|^{2-m+\epsilon}) \ \text{and}\
\beta(x)\ge-C(1+|x|)^{2-m+2\delta},
\]
for some $C,\epsilon>0$ and some $\delta<\epsilon$,
then the compact support property does not hold for $X(t)$.
\end{enumerate}
\end{theorem}

\begin{remark}
By Theorem EP1, to prove Theorem 1, it is necessary and sufficient
to show that if $\alpha$ is
as in part (1) of the theorem, then
there is no nontrivial solution
to \eqref{evoequ} with initial data $f=0$, while if $\alpha$ is as in part (2) of the theorem,
then there is such a nontrivial solution. In the case that
$L=\frac12\Delta$, and for part (2), $\beta\ge0$,
 this result was obtained in \cite[Theorem 7]{EP03}.
An alternative, more purely probabilistic proof  which does not
rely on Theorem EP1 can be found in \cite{R04} for the case $L=\frac12\Delta$
and $\beta=0$.
\end{remark}

A heuristic, qualitative understanding of Theorem \ref{smallalpha} is given at the
end of this section.

As a complement to Theorem \ref{smallalpha}, we note the
following result \cite{E00, EP03}.

\begin{EP3} Let $p\in(1,2]$.
\begin{enumerate}
\item Let $d\ge2$ and let
\[
L=A(x)\Delta,\
\text{where} \ A(x)\ge C(1+|x|)^m,\
\text{for some}\ C>0\
\text{and}\  m>2.
\]
Assume that
\[
\sup_{x\in R^d}\alpha(x)<\infty\ \text{and}\
\beta\ge0.
\]
Then the compact support
property does not hold for $X(t)$.
\item Let $d=1$ and let
\[
L=A(x)\frac{d^2}{dx^2},\
\text{where} \ A(x)\ge C(1+|x|)^m,\
\text{for some}\ C>0\
\text{and}\  m>1+p.
\]
Assume that
\[
\sup_{x\in R^d}\alpha(x)<\infty\ \text{and}\
\beta\ge0.
\]
Then the compact support
property does not hold for $X(t)$.
\item Let $d=1$ and let
\[
L=A(x)\frac{d^2}{dx^2},\
\text{where} \ A(x)\le C(1+|x|)^m,\
\text{for some}\ C>0,\
\text{and}\  m\le1+p.
\]
Assume that
\[
\inf_{x\in R^d}\alpha(x)>0\ \text{and}\
\beta\le0.
\]
Then the compact support
property holds for $X(t)$.
\end{enumerate}

\end{EP3}

\begin{remarK 1}
It follows from Theorem EP3 that if
$d=1$ and $L=(1+|x|)^m\frac{d^2}{dx^2}$, with $m\in(2,3]$, and
say $\alpha=1$ and  $\beta=0$, then the compact support property
will depend on the particular choice of $p\in(1,2]$.
\end{remarK 1}

\begin{remarK 2}
If $d=2$ and $L=(1+|x|)^m\Delta$, with $m>2$,
and say $\alpha=1$ and $\beta=0$, then by Theorem EP3,
the compact support property does not hold, yet the underlying
diffusion does not explode since it is a time-change of a recurrent
process; namely, of two-dimensional Brownian motion.

\end{remarK 2}

\begin{proof}
By Theorem EP1, it is necessary and sufficient to show that under
the conditions of parts (1) and (2),
there exists
a nontrivial solution to \eqref{evoequ}, while under
the conditions of part (3) there does not. In the case
that $\alpha=1$ (or equivalently, any positive constant) and $\beta=0$,
this follows from
\cite[Theorem 5]{EP03}. To extend this to $\alpha$ and $\beta$ as in
the statement of the theorem, one appeals to a comparison
result which we state below \cite[Proposition 4]{EP03}.

\end{proof}

\begin{compare}\label{comparison}
Assume that
$$
\beta_1\le \beta_2
$$
and
$$
0<\alpha_2\le\alpha_1.
$$
If uniqueness holds for
\eqref{evoequ} with initial data $f=0$
when $\beta=\beta_2$ and $\alpha=\alpha_2$, then
uniqueness also
holds when $\beta=\beta_1$ and $\alpha=\alpha_1$.
Thus, from Theorem EP1,
 if the compact support property holds for $\beta=\beta_2$
and $\alpha=\alpha_2$, then it also holds for
$\beta=\beta_1$ and $\alpha=\alpha_1$.
\end{compare}

Theorem \ref{smallalpha} and Theorem EP3 demonstrate the effect of
the underlying diffusion $Y(t)$ on the compact support property in
the case that $L$ is comparable to
$(1+|x|)^m\Delta$.
 We now consider more generally the effect of the
underlying diffusion process on the compact support property. We
begin with the following result which combines  \cite[Theorem
3]{EP03} with Theorem EP1.
\begin{EP4}
Let $p\in(1,2]$ and
assume that the underlying diffusion process $Y(t)$
 explodes. Assume in addition that
\[
\inf_{x\in R^d}\frac{\beta(x)}{\alpha(x)}>0.
\]
Then there exists a nontrivial solution to
\eqref{evoequ} with initial data $f=0$;
thus, by Theorem EP1,  the compact support property does not hold.
\end{EP4}

In particular, it follows from Theorem EP4 that if $Y(t)$ is explosive and  $\sup_{x\in R^d}\alpha<\infty$,
then a sufficient condition for the compact
support property to fail is that
$\inf_{x\in R^d}\beta (x)>0$.
It turns out that this result can be extended significantly.
We will prove  the following key result.

\begin{theorem}\label{equivbetas}
Let $p\in(1,2]$ and let $L$ and $\alpha$ be arbitrary.
Assume that
$$
\sup_{x\in R^d}|\beta_1-\beta_2|<\infty.
$$
Then uniqueness holds for \eqref{evoequ} with initial data $f=0$
for  $L,\alpha$ and $\beta_1$ if and only if
it holds for $L,\alpha$ and $\beta_2$. Thus, from Theorem EP1, the compact
support property holds with $\beta_1$ if and only if it holds with $\beta_2$.
\end{theorem}

\begin{remark}
Theorem \ref{equivbetas} states that a bounded change in the parameter $\beta$ cannot influence
the compact support property.
\end{remark}

As an immediate consequence of Theorem EP4 and Theorem \ref{equivbetas}, we obtain
the following result.

\begin{theorem}\label{explosionnocompact}
Let $p\in(1,2]$ and assume that the underlying diffusion process $Y(t)$ explodes.
Assume that
\begin{equation}
\sup_{x\in R^d}\alpha(x)<\infty\ \text{and}\ \inf_{x\in R^d}\beta(x)>-\infty.\label{alphabeta}
\end{equation}
Then there exists a nontrivial solution to \eqref{evoequ} with initial data $f=0$; thus
by Theorem EP1, the compact support property does not hold.
\end{theorem}

\begin{REmark 1}
Theorem \ref{explosionnocompact} shows that if the branching mechanism satisfies \eqref{alphabeta},
then the compact support property never holds if the underlying diffusion is explosive.
The converse is not true---an example was given in Remark 2 following Theorem EP3,
and another one appears in Remark 2 following Corollary \ref{hitting}.
\end{REmark 1}

\begin{REmark 2}
Starting from Theorem EP4,
a direct probabilistic proof of Theorem \ref{explosionnocompact}
can be given in the case $p=2$ via Dawson's Girsanov theorem for super-diffusion processes
\cite{D93}.
\end{REmark 2}

\begin{REmark 3}
Theorem \ref{explosionnocompact} is significant also from the point of view of pde's.
It states that under the condition \eqref{alphabeta}, nonuniqueness in the class of bounded
solutions to the linear equation \eqref{linequ} guarantees nonuniqueness for the semi-linear equation
\eqref{evoequ} with initial data $f=0$.

\end{REmark 3}

The next result shows that the restriction
$\sup_{x\in R^d}\alpha(x)<\infty$ in Theorem \ref{explosionnocompact}
is essential.

\begin{proposition}
Let $p\in(1,2]$.
Let $m\in (-\infty,\infty)$,
\[
L=(1+|x|)^m\Delta \ \text{in}\ R^d,
\]
$\beta=0$ and
\[
\alpha(x)\ge c(1+|x|)^{m-2},
\]
for some $c>0$.
Then the compact support property holds
for the measure-valued process $X(t)$.
However, if  $m>2$ and $d\ge3$,  the diffusion process $Y(t)$ explodes.
\end{proposition}

\begin{proof}
By the comparison result above, it suffices to consider
the case that  $\alpha=c(1+|x|)^{m-2}$.
Note that upon dividing by $(1+|x|)^m$,
 the stationary elliptic equation can be
written in the form
$\Delta W-\frac c{(1+|x|)^2}W^p=0$ in $R^d$.
 There is no nontrivial, nonnegative solution to this
equation (see, for
example \cite[Theorem 6]{EP03}). From this it follows that
uniqueness holds for \eqref{evoequ} with initial data $f=0$;
indeed, if uniqueness did not hold, and $u(x,t)$ were a nontrivial
solution,  then it would follow from the maximum principle that
$u(x,t)$ is increasing in $t$. Then $W(x)\equiv\lim_{t\to\infty}u(x,t)$
would constitute a nontrivial solution to the above stationary elliptic
equation. (For details, see \cite[Theorem 4-(ii)]{EP03}.) Since
uniqueness holds for \eqref{evoequ}, it follows from Theorem EP1
that the compact support property holds. For a proof that $Y(t)$
explodes if
 $m>2$ and $d\ge3$,
see, for example, \cite[Proposition 5
and Remark 1 following it]{EP03}).

\end{proof}

As mentioned in the introduction,
$p$ and $b$ play only a minor role in determining
whether or not the compact support property holds.
Theorem \ref{equivbetas} demonstrates the limited role played by $\beta$.
In remark 1 following Theorem EP3, we
 have seen a rather restricted example where $p$ can effect
the
compact support property. In the sequel we will give another example where
$p$ can effect the compact support property and also an example where $\beta$ can
effect the compact support property.
In order to accomplish this, we  first need to discuss how
the concept of a measure-valued process \it hitting\rm\ a point
can be formulated and understood in terms of the compact support property.
This last point is of independent interest.

Let $\mathcal
R_t=\text{cl}\big(\cup_{s\in[0,t]}~\text{supp}(X(s))\big)$ and let
$\mathcal R=\text{cl}\big(\cup_{s\ge0}~\text{supp}(X(s))\big)$.
 The random set $\mathcal
R$ is called the \it range \rm\ of $X=X(\cdot)$. A path of the
measure-valued process   is said to \it hit \rm a point
$x_0\in R^d$ if $x_0\in \mathcal R$. If $X(t)$  becomes extinct
with probability one, that is, $P_\mu(X(t)=0 \ \text{for all large
}\ t)=1$, or more generally, if $X(t)$ becomes locally extinct
with probability one, that is,
$P_\mu(X(t,B)=0 \ \text{for all large}\ t)=1$,
for each bounded  $B\subset R^d$,
 then $x_0\in \mathcal R$ if and only if $x_0\in
\mathcal R_t$ for sufficiently large $t$. Thus, we have:
\begin{equation}\label{hitpoint}
\begin{aligned}
&\text{If }\ X(t)\ \text{suffers local extinction with probability one,
then}\\
&P_\mu(X \ \text{hits} \ x_0)>0\ \text{if and only if there exists
a } \  t>0
\ \text{such that}\\
&P_\mu(\cup_{0\le s\le t}~ \text{supp}( X(s)) \ \text{is not
compactly embedded in}\ R^d-\{x_0\})>0.
\end{aligned}
\end{equation}

Now although we have assumed in this paper that the underlying
state space is $R^d$, everything goes through just as well on an
arbitrary domain $D\subset R^d$ \cite{EP99}. Of course now, the
compact support property is defined with respect to the domain
$D$, and  the underlying diffusion will explode if it hits
$\partial D$ in finite time.
In particular, Theorem EP1 still holds with $R^d$ replaced by $D$
\cite{EP99}.

In light of the above observations, consider a measure-valued
process $X(t)$
 corresponding to the log-Laplace
equation \eqref{evoequ} on $R^d$ with $d\ge2$. The underlying
diffusion process $Y(t)$ on $R^d$ corresponds to the operator $L$
on $R^d$. Let $\hat Y(t)$ denote the diffusion process on the
domain  $D=R^d-\{x_0\}$ with absorption at $x_0$ and
 corresponding
to the same operator $L$. (Note that if $x_0$ is polar for $Y(t)$,
then $Y(t)$ and $\hat Y(t)$
 coincide  when started from $x\neq x_0$. In fact, $x_0$ is always polar
under the assumptions we have placed on the coefficients of $L$ \cite{D98}.)
Let $\hat X(t)$
denote the measure-valued process corresponding to the
log-Laplace equation \eqref{evoequ}, but with $R^d$ replaced by
$D=R^d-\{x_0\}$. \it It  follows from \eqref{hitpoint}
that if $X(t)$ suffers local extinction with probability one, then
the measure-valued process $X(t)$
  hits the point $x_0$ with positive probability if and only
 if $\hat X(t)$ on $R^d-\{x_0\}$
 does not possess the compact support property.\rm\
Furthermore, the above discussion shows that even if $X(t)$
does not suffer local extinction with probability one, a sufficient
condition for $X(t)$  to hit the point $x_0$ with positive
probability is that $\hat X(t)$ on $R^d-\{x_0\}$ does not possess the
compact support property.

A similar analysis can be made when $d=1$. The process $\hat X(t)$ above
must be replaced by two processes, $\hat X^+(t)$ and $\hat X^-(t)$,
defined respectively on $(x_0,\infty)$ and $(-\infty,x_0)$.
The claim in italics above then holds with the requirement  on $\hat X(t)$
transferred to both $\hat X^+(t)$ and $\hat X^-(t)$. In the sequel, we will assume that $d\ge2$.

Consider now the following semi-linear equation in the punctured space $R^d-\{0\}$.
\begin{equation}\label{semi-linear}
\begin{aligned}
&u_t=\frac12\Delta u-u^p\ \text{in}\ (R^d-\{0\})\times(0,\infty);\\
& u(x,0)=0 \ \text{in}\ R^d-\{0\};\\
&u\ge0\ \text{in}\ (R^d-\{0\})\times[0,\infty).\\
\end{aligned}
\end{equation}

By Theorem EP1, the measure-valued process corresponding to
the semi-linear equation
$u_t=\frac12\Delta u-u^p$ in $R^d-\{0\}$
will possess the compact
support property if and only if \eqref{semi-linear}  has a nontrivial
solution.
The following theorem was recently proved directly in \cite{P05}; in fact it is
a particular case of a more general result in \cite{BP84}.

\begin{p}
Let $p>1$ and $d\ge2$.
\begin{enumerate}
\item If $d<\frac{2p}{p-1}$, then there exists a nontrivial solution
to \eqref{semi-linear}.
\item If $d\ge\frac{2p}{p-1}$, then there is no nontrivial
solution to \eqref{semi-linear}.
\end{enumerate}
\end{p}

\begin{remark}
In Theorem BP, $p>1$ is unrestricted, even though of course
there is no probabilistic import when $p>2$.

\end{remark}

The following result will be proved by the method used in \cite{P05}
to prove Theorem BP.

\begin{theorem}\label{beta}
Let $p>1$ and $d\ge2$.
Let $X(t)$ denote the measure-valued process
on the punctured space $R^d-\{0\}$
corresponding to the semi-linear equation
$u_t=\frac12\Delta u+\beta u-u^p\ \text{in}\ (R^d-\{0\})\times(0,\infty)$.
Consider the Cauchy problem
\begin{equation}\label{singular}
\begin{aligned}
&u_t=\frac12\Delta u+\beta u-u^p\ \text{in}\ (R^d-\{0\})\times(0,\infty);\\
& u(x,0)=0 \ \text{in}\ R^d-\{0\};\\
&u\ge0\ \text{in}\ (R^d-\{0\})\times[0,\infty).\\
\end{aligned}
\end{equation}
Assume that
\[
d<\frac{2p}{p-1}.
\]
Let
\[
\beta_0=\frac{d(p-1)-2p}{(p-1)^2}<0.
\]
\begin{enumerate}
\item If
\[
\beta(x)\ge\frac{\beta_0+\kappa}{|x|^2},\ \text{for some}\ \kappa\in(0,-\beta_0],
\]
then there exists a nontrivial solution to \eqref{singular}; hence, the compact
support property does not hold for $X(t)$.
\item
If
\[
\limsup_{x\to0}|x|^2\beta(x)<\beta_0,
\]
then there is no nontrivial
solution to \eqref{singular}; hence, the compact support property
holds for $X(t)$.

\end{enumerate}

\end{theorem}

\begin{remark}
The restriction $\kappa\le-\beta_0$ is made to ensure that $\beta$ is bounded from
above.
\end{remark}

We have the following corollary of Theorem BP and  Theorem \ref{beta}.

\begin{corollary}\label{hitting}
Let $X(t)$ denote the measure-valued process on all of $R^d$, $d\ge2$,
corresponding to the semi-linear equation
\[
u_t=\frac12\Delta u+\beta u-u^p \ \text{in}\ R^d\times(0,\infty).
\]
\begin{enumerate}
\item
If $\beta$ is bounded from below and $d<\frac{2p}{p-1}$, then $X(t)$ hits any point $x_0$
with positive probability;

\item
If  $\beta\le0$ and $d\ge\frac{2p}{p-1}$, then $X(t)$  hits any point
$x_0$ with probability 0;

\item If $\beta\le 0$, $d<\frac{2p}{p-1}$
and $\beta$ has a singularity
at the origin such that
\[
\limsup_{x\to0}|x|^2\beta(x)<\frac{d(p-1)-2p}{(p-1)^2},
\]
then $X(t)$ hits 0 with probability 0.

\end{enumerate}

\end{corollary}

\begin{proof}
When $\beta=0$, part (1) follows immediately from Theorem BP,
Theorem EP1 and the discussion preceding Theorem BP. For the general case
one appeals to Theorem \ref{equivbetas}, which holds just as well
for a punctured space. (Recall that we are always assuming in this paper
that $\beta$ is bounded.)

When $\beta=0$, it is well-known that the super-Brownian motion in the statement
of the corollary suffers extinction with probability one. By comparison, this
also holds when $\beta\le0$ \cite{EP99}. Thus, part (2) follows from
Theorem BP, Theorem EP1 and the discussion preceding Theorem BP, while part (3)
follows from Theorem \ref{beta}, Theorem EP1 and the discussion preceding Theorem BP.
\end{proof}

\begin{remark 1} When $\beta=0$,
the results in parts (1) and (2) of Corollary \ref{hitting}
state that critical, super-Brownian motion hits a point with positive
probability if $d<\frac{2p}{p-1}$, and with zero
probability if $d\ge\frac{2p}{p-1}$. This result
can be found in \cite{DIP89} for the case $p=2$.
For $p\in(1,2]$, it can be found in
\cite{D91} or in
\cite{DK96},
which exploit the method of removable singularities
for \it elliptic\rm\ equations, developed in \cite{BV80} and \cite{V81}.
The approach here is via the parabolic equation.

\end{remark 1}

\begin{remark 2}
Theorem BP gives an example where the measure-valued process does not
possess the compact support property even though the underlying diffusion
process does not explode. Indeed, the underlying diffusion is Brownian
motion in $R^d-\{0\}$. Since singletons are polar for multi-dimensional
Brownian motion, this process does not explode. However for $d<\frac{2p}{p-1}$,
the compact support property fails for the measure-valued process.

\end{remark 2}

\begin{remark 3}
In order to obtain an example of the phenomenon occurring in Remark 2
when the state space is the whole space,
and in order to see how the parameter $p$ and the mass creation
parameter $\beta$ can affect the compact support property when
the state space is the whole space, we convert the set-up in
Theorem \ref{beta} and Theorem BP
to the state space $R=(-\infty,\infty)$
by considering just the radial variable, $r=|x|$, and then
 making
a change of variables, say, $z=\frac1r-r$.
One obtains an operator $L=\frac12a(z)\frac{d^2}{dz^2}+b(z)\frac d{dz}$
on $R$,
where
\[
\begin{aligned}&a(z)\sim z^4 \ \text{as}\ \ z\to\infty, \
\lim_{z\to-\infty}a(z)=1;\\
& b(z)\sim \frac{3-d}2z^3,\  \text{as}\ \ z\to\infty, \
\lim_{z\to-\infty}b(z)=0.
\end{aligned}
\]
Also, one has $\alpha=1$. One has $\beta=0$ in Theorem BP and one has
\begin{subequations}\label{joint}
\begin{equation}\label{part1}
\beta(z)\ge\frac{4(\beta_0+\kappa)}{\big((z^2+4)^\frac12-z)\big)^2}
\ \ \text{in part (1) of Theorem \ref{beta}},
\end{equation}
\begin{equation}\label{part2}
\limsup_{z\to\infty}\frac14\big((z^2+4)^\frac12-z)\big)^2\beta(z)
<\beta_0 \ \ \text{in part (2)
of Theorem \ref{beta}}.
\end{equation}

\end{subequations}
\noindent It follows then that $\liminf_{z\to\infty}\frac{\beta(z)}{z^2}\ge\beta_0+\kappa$
in part (1) of Theorem \ref{beta} and $\limsup_{z\to\infty}\frac{\beta(z)}{z^2}<\beta_0$
in part (2) of Theorem \ref{beta}.

Consider first the phenomenon mentioned in Remark 2
with regard to Theorem BP.
After the change of variables, the state space is $R$ and the operator $L$ on $R$ depends
on the parameter $d$. The one-dimensional diffusion corresponding to $L$ is nonexplosive for
all $d\ge2$, as it inherits this property from the original process before
the change of variables. Also $\alpha=1$ and $\beta=0$.
However, if $d<\frac{2p}{p-1}$, then the compact support property fails.
(Another example of this phenomenon was presented in Remark 2 following
Theorem EP3.)

With the same setup as in the previous paragraph, we also see how the parameter
$p$ affects the compact support property---the property will hold if and only
if $p\ge\frac d{d-2}$.

Now consider the above change of variables applied to Theorem \ref{beta}. We see
that if $d<\frac{2p}{p-1}$, then $\beta$ affects the compact support
property---it holds if $\beta$ satisfies \eqref{part2} and does
not hold if $\beta$ satisfies \eqref{part1}.
As Theorem \ref{equivbetas} guarantees, an unbounded change in $\beta$ was needed to effect
a change in the compact support property.

\end{remark 3}

\medskip

Theorems \ref{smallalpha}, \ref{equivbetas} and \ref{beta}
will
be proved successively in the  sections that follow.

\medskip

We now turn to an intuitive probabilistic understanding of
the role of the branching in
Theorem 1.
We recall the particle
process approximation to the measure-valued process
in the case that $\alpha$ and $\beta$ are bounded.
Consider first the case $p=2$, the case
in which the offspring distribution has finite variance.
Let  $ \hat R^d=R^d\cup \{\Delta\}$
denote the one-point compactification of $R^d$. One may consider
the diffusion process $Y(t)$ to live on $\hat R^d$; if $Y(t)$ does not
explode,  then it never reaches $\Delta$, while if $Y(t)$
does  explode, then it enters the state $\Delta$ upon leaving $R^d$,
and remains there forever.
For each positive integer n, consider $N_{n}$
particles, each of mass $\frac{1}n$, starting at points $y_{i}^{(n)}(0)\in
 R^d,i=1,2,\dots,N_{n},$ and performing independent branching diffusion according
to the  process $Y(t)$, with branching rate $cn,c>0$, and spatially
dependent branching
distribution $\{p_{k}^{(n)}(y)\}_{k=0}^{\infty}$, where

\[
\sum_{k=0}^{\infty}kp_{k}^{(n)}(y)=1+\frac{\gamma(y)}n+o(\frac1n),\
\text{as}\ n\to\infty;
\]
\[
\sum_{k=0}^{\infty}(k-1)^{2}p_{k}^{(n)}(y)=m(y)+o(1),\ \text{as}\ n\to
\infty, \ \text{uniformly in}\ y,
\]
with $m,\gamma \in C^\alpha( R^d)$ and $m(y)>0$. Let $N_{n}(t)$ denote the number
of particles alive at time
t and denote their positions by $\{Y_{i}^{(n)}(t)\}_{i=1}^{N_{n}(t)}$. Denote
by $\mathcal{M}_{F}( R^d)$ ($\mathcal{M}_{F}(\hat R^d)$) the space of finite measures
on $ R^d$ ($\hat R^d$). Define an $\mathcal{M}_{F}(\hat R^d)$- valued process $X_{n}(t)$
by $X_{n}(t)=\frac{1}n\sum_{1}^{N_{n}(t)}\delta_{Y_{i}^{(n)}(t)}(\cdot)$.
Note that $X_n$ is c\`adl\`ag.
Denote by $P^{(n)}$ the probability measure
corresponding to $\{X_{n}(t), 0\le t<\infty\}$
on $D([0,\infty),\mathcal M_F(\hat R^d))$, the space of c\`adl\`ag paths
with the Skorohod topology,

Assume that $m(y)$ and $\gamma(y)$ are  bounded from above.
One can show
that if $w-\lim_{n\to\infty} X_n(0)=\mu\in \mathcal M_F( R^d)$,
then
$P^*_\mu=w-\lim_{n\to\infty}P^{(n)}$ exists
in  $D([0,\infty),\mathcal M_F(\hat R^d))$.
Furthermore, the measure
$P_\mu^*$ restricted to $D([0,\infty),\mathcal M_F(R^d))$
satisfies
 \eqref{Def:MVD} and \eqref{evoequ}
with
$\beta(y)=c\gamma(y)$, $\alpha(y)=\frac{1}{2}cm(y)$
and $p=2$ (see \cite{EP99}). Denoting this restriction by $P_\mu^*|_{R^d}$,
it then follows that  $P_\mu=P_\mu^*|_{R^d}$. In fact, one
can show that
$P_\mu$ is supported on the space of continuous
paths, $ C([0,\infty),\mathcal M_F( R^d))$.

One should think of $\beta$ and $\alpha$ as the \it mass creation\rm\
and  the \it variance\rm\ parameters  respectively of the branching.

For the case that $p\in(1,2)$, one cooks up a
sequence of  distributions $\{p_k^{(n)}(y)\}_{n=1}^\infty$ for which
the generating functions $\Phi^{(n)}(s;y)=\sum_{k=0}^\infty
p_k^{(n)}(y)s^y$ satisfy
$\lim_{n\to\infty}n^p\big(\Phi^{(n)}(1-\frac\lambda n;y)-(1-\frac\lambda n)
\big)=\alpha(y)\lambda^p-\beta(y)\lambda$.
These offspring
 distributions $\{p_k^{(n)}(y)\}_{n=1}^\infty$
will possess all moments smaller than $p$.
Again, $\beta$ can be thought of as the mass creation parameter.
As above, each particle is given mass $\frac1n$, but in the present case,
the branching rate is $n^{p-1}$.
The same construction and conclusion as above holds, although
in this case the paths are not continuous, but only c\`adl\`ag
\cite{D93}.

With the above set-up, we can now give some intuition concerning Theorem 1.
Consider two particular cases of the above construction.
In both cases we will assume that
at time 0 there are $n$ particles, all positioned at $y=0$; that is,
$N_n=n$ and $y_i^{(n)}=0$, for $i=1,2,\ldots n$. Then  the
initial measure, both for the approximating process and the limiting
one, will be $\mu=\delta_0$.
We will also assume that the diffusion $Y(t)$ does not explode.

The first case is the  completely trivial case in which there is
\it no branching
at all.\rm\ This  degenerate case corresponds to $\beta=\alpha=0$
(and thus does not actually fit into the above set-up).
In this case, $X_n(t)$ is  a random  probability measure with
$n$ atoms of mass $\frac1n$ positioned at $n$ IID points, distributed
according to the distribution  of $Y(t)$. Thus, by the law of large
numbers, $X(t)=w-\lim_{n\to\infty}X_n(t)$ is the deterministic measure
$\text{dist}(Y(t))$.
Since $\text{dist}(Y(t))$ is not compactly supported for $t>0$, it follows that the
compact support property does not hold for $X(t)$ in this trivial case.

Now consider
the case of critical, binary branching; that is,  $p_0^{(n)}=p_2^{(n)}=\frac12$.
Letting $c=1$, it then follows that  $\beta=0$ and $\alpha=\frac12$.
In this case, it is well-known that
for any $s<1$, if one
lets $p_n(s,M)$ denote
 the probability that
all the mass at time 1 in the
approximate measure-valued process $X_n(\cdot)$
descends from no more than $M$ ancestors
alive at time $s$, then $\lim_{M\to\infty}\lim_{n\to\infty}p_n(s,M)=1$.
Thus, since all the particles  alive at time 1 are coming
from a finite number of ancestors at any time $s$, these particles
are correlated, and the law of large numbers does not apply, allowing
for compact support property to hold.

The above discussion suggests that one way for the compact support
property to break down is for the branching mechanism
to be spatially dependent and to decay sufficiently fast as $|x|\to\infty$
so that the law of large numbers will come into play.
Furthermore, the faster the diffusion is, the more quickly individual
particles that begin together become statistically uncorrelated,
so one might expect that the  stronger the
diffusion, the weaker the threshold  on the
decay rate in order for the compact support property to break down.
If the diffusion process $Y(t)$ corresponds to the operator
$L=\frac12(1+|x|)^m\Delta$, for $m\in[0,2]$,
then for $m=0$ one obtains Brownian motion, while for $m\in(0,2]$, one
obtains a time-changed Brownian motion with the diffusion sped up,
the speed increasing in $m$.
Theorem one shows that for such an underlying diffusion, for any
$\epsilon>0$,
the rate $\exp(-|x|^{2-m+\epsilon})$ is sufficiently
fast to cause the compact support property to fail, but  the rate
$\exp(-|x|^{2-m})$ is not fast enough.

\section{Proof of Theorem \ref{smallalpha}}\label{S:smallalpha}
To make the calculations simpler, we will prove
the theorem in the case that $L=A(x)\Delta$, where
$C_0^{-1}(1+|x|)^m\le A(x)\le C_0(1+|x|)^m$, for some $m\in[0,2]$.
The general case follows in the same fashion.
\begin{proof}[Proof of part (1)]
 Let $u(x,t)$ be any solution of \eqref{evoequ} with initial data
$g=0$.
By Theorem EP1, we need to show that $u\equiv0$.
Define $\hat U(x,t)$ through the equality
$u(x,t)=\hat U(x,t)\exp(\lambda(1+|x|^2)^{\frac{2-m}2}(t+\delta))$,
for some $\lambda,\delta>0$.
Then
\begin{equation}\label{Laplacian}
\begin{aligned}
&\exp(-\lambda(1+|x|^2)^{\frac{2-m}2}(t+\delta))A(x)\Delta u=\\
&A(x)\Big(\Delta \hat U+2\lambda(2-m)(1+|x|^2)^{-\frac m2}(t+\delta) x\cdot\nabla
\hat U\Big)\\
&+A(x)\Big((t+\delta)^2\lambda^2(2-m)^2(1+|x|^2)^{-m}|x|^2\\
&+(t+\delta)\lambda d(2-m)(1+|x|^2)^{-\frac m2}-(t+\delta)\lambda(2-m)m(1+|x|^2)^{-\frac m2-1}
|x|^2)\Big)\hat U
\end{aligned}
\end{equation}
Also, using the bound on $\alpha$ in the statement of the theorem, we have
\begin{equation}\label{otherterms}
\begin{aligned}
&\exp(-\lambda(1+|x|^2)^{\frac{2-m}2}(t+\delta))\Big(\beta u-\alpha u^p-u_t\Big)=\\
&\Big(\beta-\lambda(1+|x|^2)^{\frac {2-m}2}\Big)\hat U-
\alpha\exp\big((p-1)\lambda(1+|x|^2)^{\frac{2-m}2}(t+\delta)\big)
\hat U^p.
\end{aligned}
\end{equation}

Since $u$ is  a solution to \eqref{evoequ}, the sum of the left hand sides
of \eqref{Laplacian} and \eqref{otherterms} is equal to 0.
Consider now the sum of the terms on the right hand sides of \eqref{Laplacian}
and \eqref{otherterms} with the variable $t$ restricted by
$0\le t\le \delta$.
By assumption, $\alpha(x)\ge C_1\exp(-C_2|x|^{2-m})$. Thus,
the coefficient of $\hat U^p$ will be bounded away
from 0 if
\begin{equation}\label{lambdadelta1}
\lambda\delta=\frac{C_2}{p-1}.
\end{equation}
In the case that $m=2$, the coefficient of $\hat U$ is bounded from above. Otherwise,
the two unbounded terms in the coefficient of $\hat U$
are $A(x)(t+\delta)^2\lambda^2(2-m)^2(1+|x|^2)^{-m}|x|^2$
and $-\lambda(1+|x|^2)^{\frac{2-m}2}$.
By assumption, there exists a $C_4>0$ such that $A(x)\le C_4(1+|x|^2)^\frac m2$.
Thus, in order to guarantee that
the coefficient of $\hat U$ is bounded from above, it suffices to have
$\lambda= C_4(2\delta)^2\lambda^2(2-m)^2$. Using \eqref{lambdadelta1} to
substitute for $(\lambda\delta)^2$ on the right hand side  above, we have
\begin{equation}\label{lambdadelta2}
\lambda=4C_4(2-m)^2C_2^2(p-1)^{-2}.
\end{equation}
With $\lambda$ and $\delta$
chosen as in \eqref{lambdadelta1} and \eqref{lambdadelta2}, it then follows that
\begin{equation}\label{subsol}
\begin{aligned}
&A(x)\Big(\Delta \hat U+2\lambda(2-m)(1+|x|^2)^{-\frac m2}(t+\delta) x\cdot\nabla
\hat U\Big)\\
&+c_1\hat U-c_2\hat U^p\ge0,\ \text{for}\ (x,t)\in R^n\times[0,\delta],
\ \text{for some} \ c_1,c_2>0.
\end{aligned}
\end{equation}

Let
$M_{R,K}(x,t)=(1+|x|)^{\frac2{p-1}}
(R-|x|)^{-\frac2{p-1}}\exp(K(t+1))$, for $(x,t)\in B_R\times(0,\infty)$.
In \cite[proof of Theorem 2]{EP03}, it was shown that
for any operator $\mathcal L$
satisfying the conditions of Theorem EP2,
there exists a $K>0$ such that for all
$R>0$
\begin{equation}\label{testfunction}
\mathcal LM_{R,K}(x,t)\le0,
\ (x,t)\in B_R\times(0,\infty),
\end{equation}
where $B_R$ denotes the ball of radius $R$ centered at the origin.
In particular the operator
$\mathcal A=
A(x)\Big(\Delta  +2\lambda(2-m)(1+|x|^2)^{-\frac m2}(t+\delta) x\cdot\nabla
\Big)$ satisfies the conditions of Theorem EP2,
except for the fact that it is time inhomogeneous.
The time inhomogeneity causes no problem since we are only considering
$t\in[0,\delta]$ and since for fixed $x$, everything in sight is uniformly
bounded for $t\in[0,\delta]$.
Thus the proof
of \eqref{testfunction} shows that
\begin{equation}\label{supersol}
\mathcal AM_{R,K}(x,t)\le0,
\ (x,t)\in B_R\times(0,\delta].
\end{equation}
Since $0=\hat U(x,0)\le M_{R,K}(x,0)$ and
$\hat U(y,t)\le M_{R,K}(y,t)=\infty$, for $y\in \partial B_R$,
it follows from \eqref{subsol}, \eqref{supersol}
and the maximum principle for semi-linear equations \cite[Proposition 1]{EP03}
that
\begin{equation}
\hat U(x,t)\le M_{R,K}(x,t), \text{for}\ (x,t)\in B_R\times[0,\delta].
\end{equation}
Letting $R\to\infty$, we conclude that $\hat U(x,t)=0$, for $(x,t)\in R^d
\times(0,\delta]$.
Thus, we also have $u(x,t)=0$, for $(x,t)\in R^d\times[0,\delta]$.
Since $u$ satisfies a time homogeneous equation, we conclude
that in fact $u(x,t)=0$, for $(x,t)\in R^d\times[0,\infty)$.
This completes the proof of part (1).
\medskip

\noindent  \it Proof of part (2).\rm\
By Theorem EP1, we need to show that there is a nontrivial solution
to \eqref{evoequ} with initial data $g=0$.
Let $u(x,t)$ be any solution of \eqref{evoequ} with initial data
$g=0$. Define $\hat U(x,t)$ through the equality
$u(x,t)=\hat U(x,t)\exp(\lambda(1+|x|^2)^{\frac{2-m+\kappa}2})$,
where $\kappa\in(\delta,\epsilon)$, and $\epsilon$ and $\delta$
are as in the statement of the theorem.
Now \eqref{Laplacian} and \eqref{otherterms}
hold with the following changes: (i) $m$ is replaced by $m-
\kappa$; (ii) $(t+\delta)$ is replaced by 1;
(iii) the term $\lambda(1+|x|^2)^{\frac{2-m}2}$ in \eqref{otherterms}
is deleted. As before,
 the sum of the left hand sides
of \eqref{Laplacian} and \eqref{otherterms} is equal to 0.
By assumption, $\alpha(x)\le C\exp(-|x|^{2-m+\epsilon})$.
Thus, the coefficient of $\hat U^p$ in the amended version
of \eqref{otherterms} is bounded from above, since $\kappa<\epsilon$.
The coefficient of $\hat U$ in the amended version of
\eqref{otherterms} is $\beta$, and by assumption, there exists a
$C_5>0$ such that
$\beta\ge-C_5(1+|x|^2)^{\frac{2-m+2\delta}2}$.
The coefficient of $\hat U$ in the amended version of \eqref{Laplacian}
is $A(x)\Big(\lambda^2(2-m+\kappa)^2(1+|x|^2)^{-m+\kappa}|x|^2
+\lambda d(2-m+\kappa)(1+|x|^2)^{-\frac {m-\kappa}2}-
\lambda(2-m+\kappa)(m-\kappa)(1+|x|^2)^{-\frac {m-\kappa}2-1}
|x|^2)\Big)$.
By assumption, there exists a $C_3>0$
such that $A(x)\ge C_3(1+|x|^2)^{\frac m2}$.
It is easy to check that
 by choosing  $\lambda$ sufficiently large,
the factor multiplying $A(x)$ above  will be  bounded from below
by $\frac{2C_5}{C_3}(1+|x|^2)^{\frac{2-2m+2\kappa}2}$.
(For  $|x|\le\frac12$, use the second term in the parentheses, and for
$|x|>\frac12$ use the first and third terms.)
Thus, the coefficient of $\hat U$ in
the amended version of \eqref{Laplacian} is
greater or equal to $2C_5(1+|x|^2)^{\frac{2-m+2\kappa}2}$.
Since $\kappa>\delta$, it follows that
the coefficient of $\hat U$ from the sum of the right hand sides
of \eqref{Laplacian} and \eqref{otherterms} is bounded below by
a positive constant.

The above analysis shows that $\hat U$ satisfies the equation
\begin{equation}\label{hatU}
\begin{aligned}
&\hat U_t=A(x)
\Big(\Delta \hat U+2\lambda(2-m+\kappa)(1+|x|^2)^{-\frac {m-\kappa}2} x\cdot\nabla
\hat U\Big)+\hat \beta\hat U-\hat \alpha \hat U^p;\\
&\hat U(x,0)=0,
\end{aligned}
\end{equation}
where $\hat \beta\ge C_6>0$ and $\hat \alpha\le C_7$,
and that uniqueness for the original equation is equivalent to uniqueness
for \eqref{hatU}.
We will show below that the diffusion process corresponding to the operator
$A(x)
\Big(\Delta +2\lambda(2-m+\kappa)(1+|x|^2)^{-\frac {m-\kappa}2} x\cdot\nabla
 \Big)$ explodes. Thus, it follows from Theorem EP4 that
 uniqueness does not hold for
\eqref{hatU}. Consequently, uniqueness does not hold for the original equation
with initial data $g=0$; thus, by Theorem EP1, the compact support property
does not hold.

It remains to show that the diffusion
corresponding to the operator \linebreak
$A(x)
\Big(\Delta +2\lambda(2-m+\kappa)(1+|x|^2)^{-\frac {m-\kappa}2} x\cdot\nabla
 \Big)$
explodes. The diffusion in question is a time change of the diffusion
corresponding to
$\Delta+2\lambda(2-m+\kappa)(1+|x|^2)^{-\frac {m-\kappa}2} x\cdot\nabla$.
It is not hard to show that since
 $A(x)\ge C_3(1+|x|^2)^{\frac m2}$, explosion will occur
 for the diffusion in question if it occurs for the diffusion corresponding
 to the operator
 $C_3(1+|x|^2)^{\frac m2}
\Big(\Delta +2\lambda(2-m+\kappa)(1+|x|^2)^{-\frac {m-\kappa}2} x\cdot\nabla
 \Big)$. This latter operator is radially symmetric, and its radial
 component is of the form
 $p(r)\frac{d^2}{dr^2}+q(r)\frac d{dr}$, with
 $p(r)$ satisfying $c_1(1+r)^m\le p(r)\le c_2(1+r)^m$, for all $r\ge0$,
 and $q(r)$ satisfying $c_3(1+r)^{1+\kappa}\le q(r)\le c_4(1+r)^{1+\kappa}$,
 for $r>1$, where $c_1, c_2,c_3,c_4>0$. This diffusion explodes by the Feller criterion
 \cite[Theorem 5.1.5]{P95}.

\end{proof}

\section{Proof of Theorem \ref{equivbetas}}

\begin{proof}
By assumption, there exists a $B>0$ such that $\beta_2\le \beta_1+B$.
For $i=1,2$, consider the following parabolic equations:
\begin{equation}\label{beta_i}
\begin{aligned}&
u_t=Lu+\beta_i u-\alpha u^{p}\ \text{on}\ R^d\times
(0,\infty);\\
&u(x,0)= 0,\  x\in R^d;\\
&u\ge0.
\end{aligned}
\end{equation}
To prove the theorem it suffices to show that if \eqref{beta_i} has a non-zero
solution when $\beta_i=\beta_2$, then it also has a non-zero solution when $\beta_i=\beta_1$.

Suppose to the contrary that
a non-zero solution to \eqref{beta_i} exists when $\beta_i=\beta_2$ but not when $\beta_i=\beta_1$.
Let $B_m$ denote the ball of radius $m$ centered at the origin in $R^d$.
For $i=1,2$, consider the functions $u^{(i)}_m, m=1,2,...$
(constructed in \cite{EP99} for $p=2$; see also \cite{EP03}
for $p\in(1,2]$),
where $u^{(i)}_m$ solves the equation
\begin{equation}\label{approx}
\begin{aligned}
&\frac{\partial u_{m}}{\partial t}=Lu_m+\beta_i u-\alpha u_m^{p}\ \ \text{on}\ B_m
\times(0,\infty):\\
&\lim_{x\to\partial B_m}u_m(x,t)=\infty,\ t>0;\\
&u_m(x,0)=0,\ x\in R^d.\\
&u_m\ge 0.
\end{aligned}
\end{equation}
Since we have assumed that when $\beta_i=\beta_1$, the only solution to
\eqref{beta_i} is the zero function, it follows from
the construction in \cite{EP99}, \cite{EP03},
that
\begin{equation}\label{zerolimit}
\lim_{m\to\infty}u_m^{(1)}=0.
\end{equation}
The same construction shows that since we have assumed that there exists a non-zero solution
to \eqref{beta_i} when $\beta_i=\beta_2$, we have
\begin{equation}\label{nonzerolimit}
\lim_{m\to\infty}u_m^{(2)}\neq0.
\end{equation}

Define $v_m(x,t)=e^{B t}u^{(1)}_m(x,t) \ \text{on}\  B_m
\times(0,\infty)$.
Using the fact  that
$B>0$ and $ p>1$  along with \eqref{approx}
 gives
\begin{equation}\label{supersolution}
Lv_m+\beta_2  v_m-v_m^p-\frac{\partial v_m}{\partial t}\le Lv_m+(\beta_1+B)
v_m-v_m^p-\frac{\partial v_m}{\partial t}\le 0.
\end{equation}
We have the boundary and initial conditions
\begin{equation}
\begin{aligned}
&\lim_{x\to\partial B_m}v_m(x,t)=\infty,\ t>0;\\
&v_m(x,0)=0,\ x\in B_m.
\end{aligned}
\end{equation}
Then an application of the semi-linear parabolic maximum principle
(\cite[Proposition
1]{EP03}) gives
\begin{equation}\label{smaller}
v_m\ge u_m^{(2)}.
\end{equation}
(In fact, one has compare $v_m$ to $u^{(2)}_{n,m}$,
 where $u^{(2)}_{n,m}$
is the minimal nonnegative solution to the inhomogeneous
semi-linear equation
$$
\begin{aligned}
&\frac{\partial u^{(2)}_{n,m}}{\partial t}=Lu^{(2)}_{n,m}+\beta_2
u^{(2)}_{n,m}-\alpha (u^{(2)}_{n,m})^{p}- \psi_{n,m}\ \text{on}\ \
B_{2m}\times(0,\infty);\\
&u^{(2)}_{n,m}(x,0)=0, \ x\in B_{2m};\\
&u^{(2)}_{n,m}\ge 0,
\end{aligned}
$$
and $0\le \psi_{n,m}\le n$ is a  function on $R^d$ vanishing on $B_m$
and equal to $n$ on $\{x\in R^d\mid \mathrm{dist}(x,B_m)\ge
1/n\}$. One has
$u_m^{(2)}=\lim_{n\to\infty}u^{(2)}_{n,m}$. See the proof of  \cite[Theorem 3.4]{EP99} for more
elaboration; see also \cite{EP03}.)

By \eqref{zerolimit} and  the definition of $v_m$, we have
$$
\lim_{m\to\infty}v_m=0;
$$
thus, by (\ref{smaller}) we have
$$
\lim_{m\to\infty}u^{(2)}_m=0,
$$
which contradicts \eqref{nonzerolimit}.

\end{proof}

\section{Proof of Theorem \ref{beta}}

\begin{proof}[Proof of part (1).]
We will show that uniqueness does not hold for
\eqref{singular}. Then by Theorem EP1, the compact support property does not hold.
By the comparison result stated after the proof of Theorem EP3,
we may assume that $\beta=\frac{\beta_0+\kappa}{|x|^2}$.
Since the problem
is now radially symmetric, it suffices to show that there exists a nontrivial
solution to the radially symmetric equation
\begin{equation}\label{onedim}
\begin{aligned}
&u_t=\frac12u_{rr}+\frac {d-1}{2r}u_r+\frac{\beta_0+k}{r^2}u-u^p,\ r\in(0,\infty), \ t>0;\\
&
u(r,0)=0, \ r\in(0,\infty);\\
&u\ge0,\   r\in(0,\infty), \ t\ge0.
\end{aligned}
\end{equation}
The function $W(x)=\kappa^{\frac1{p-1}}r^{-\frac2{p-1}}$
is a positive, stationary solution
of the parabolic equation
$u_t=\frac12u_{rr}+\frac {d-1}{2r}u_r+\frac{\beta_0+\kappa}{r^2}u-u^p$ in $(0,\infty)$.
By \cite[Theorem 2-ii]{P05},
the fact that there exists a nontrivial positive, stationary solution
guarantees that uniqueness does not hold for the corresponding
parabolic equation with initial data 0; that is, uniqueness does not
hold for \eqref{onedim}.
Actually, the result in
\cite{P05} is for equations with domain $R^d$, $d\ge1$, whereas
the domain here is $(0,\infty)$. One can check that the proof also
holds in a half space, but more simply, one can make the change of
variables $z=\frac1x-x$, which converts the problem to all of $R$.

\noindent \it Proof of part (2).\rm\
We will show that uniqueness holds for \eqref{singular}. Then by Theorem EP1, the
compact support property holds.
For $\epsilon$ and $R$ satisfying
$0<\epsilon<1$ and $R>1$, and for some $l\in(0,1]$, define
\begin{equation}\label{phi}
\phi_{R,\epsilon}(x)=((|x|-\epsilon)(R-|x|))^{-\frac2{p-1}}(1+|x|)^\frac2{p-1}
(1+\frac{\epsilon^l}{ |x|^l} R^\frac2{p-1}).
\end{equation}
Also, for $R$ and $\epsilon$ as above, and some $\gamma>0$, define
\begin{equation}\label{psi}
\psi_{R,\epsilon}(x,t)=\phi_{R,\epsilon}(x)\exp(\gamma(t+1)).
\end{equation}
Note that $\psi_{R,\epsilon}(x,0)>0$, for
$|x|\in(\epsilon, R)$, and  $\psi_{R,\epsilon}(x,t)
=\infty$, for $|x|=\epsilon$ and $|x|=R$.
We will show that for
all sufficiently large $R$ and all sufficiently  small $\epsilon$, and for
$\gamma$ sufficiently large and $l$ sufficiently small, independent of those
 $R$ and $\epsilon$, one has
\begin{equation}\label{supersolution}
\frac12\Delta\psi_{R,\epsilon}
+\beta \psi_{R,\epsilon}-\psi^p_{R,\epsilon}
-(\psi_{R,\epsilon})_t\le0,\ \ \text{for} \
\epsilon<|x|<R\ \ \text{and}\ t>0.
\end{equation}
It then follows from the maximum principle for semi-linear equations
 \cite[Proposition 1]{EP03}
that  every  solution $u(x,t)$ to \eqref{singular} satisfies
\begin{equation}\label{upsi}
u(x,t)\le\psi_{R,\epsilon}(x,t), \ \ \text{for}\ \epsilon<|x|<R
\ \ \text{and} \ t\in[0,\infty).
\end{equation}
Substituting \eqref{phi} and \eqref{psi} in \eqref{upsi}, letting
$\epsilon\to0$, and then letting $R\to\infty$, we conclude  that
$u(x,t)\equiv0$.
Thus, it remains to show \eqref{supersolution}.

From now on we will use radial coordinates, writing $\phi(r)$
for $\phi(x)$ with $|x|=r$ and similarly for $\psi$.
We have
\begin{equation}\label{first}
\begin{aligned}
&\exp(-\gamma(t+1))(\psi_{R,\epsilon})_r=\\
&-(\frac2{p-1})((r-\epsilon)(R-r))^{-\frac2{p-1}-1}
(R+\epsilon-2r)(1+r)^\frac2{p-1}(1+\frac{\epsilon^l}{r^l}R^\frac2{p-1})\\
&+(\frac2{p-1})((r-\epsilon)(R-r))^{-\frac2{p-1}}(1+r)^{\frac2{p-1}-1}
(1+\frac{\epsilon^l}{r^l}R^\frac2{p-1})\\
&-l((r-\epsilon)(R-r))^{-\frac2{p-1}}(1+r)^{\frac2{p-1}}
\frac{\epsilon^l}{r^{l+1}}R^\frac2{p-1},
\end{aligned}
\end{equation}
and
\begin{equation}\label{second}
\begin{aligned}
&\exp(-\gamma(t+1))\big((r-\epsilon)(R-r)\big)^{-\frac2{p-1}-2}
(\frac12\psi_{R,\epsilon})_{rr}=\\
&(\frac1{p-1})(\frac2{p-1}+1)(R+\epsilon-2r)^2(1+r)^\frac2{p-1}
(1+\frac{\epsilon^l}{r^l}R^\frac2{p-1})\\
&+(\frac2{p-1})(r-\epsilon)
(R-r)(1+r)^{\frac2{p-1}}(1+\frac{\epsilon^l}{r^l}R^\frac2{p-1})\\
&-(\frac2{p-1})^2(r-\epsilon)(R-r)(R+\epsilon-2r)(1+r)^{\frac2{p-1}-1}
(1+\frac{\epsilon^l}{r^l}R^\frac2{p-1})\\
&+l(\frac2{p-1})(r-\epsilon)(R-r)(R+\epsilon-2r)(1+r)^{\frac2{p-1}}
\frac{\epsilon^l}{r^{l+1}}R^{\frac2{p-1}}\\
&+(\frac1{p-1})(\frac2{p-1}-1)((r-\epsilon)(R-r))^2
(1+r)^{\frac2{p-1}-2}
(1+\frac{\epsilon^l}{r^l}R^\frac2{p-1})\\
&-l(\frac2{p-1})((r-\epsilon)(R-r))^2
(1+r)^{\frac2{p-1}-1}
\frac{\epsilon^l} {r^{l+1}}R^\frac2{p-1}\\
&+\frac{l(l+1)}2((r-\epsilon)(R-r))^2(1+r)^\frac2{p-1}
\frac\epsilon{r^{l+2}}R^\frac2{p-1}.
\end{aligned}
\end{equation}
Using \eqref{psi}, \eqref{first} and the fact that $\frac2{p-1}+2=\frac{2p}{p-1}$,
we have
\begin{equation}\label{other}
\begin{aligned}
&\exp(-\gamma(t+1))\big((r-\epsilon)(R-r)\big)^{-\frac2{p-1}-2}\times\\
&\Big(\frac {d-1}{2r}(\psi_{R,\epsilon})_r+\beta\psi_{R,\epsilon}-\psi_{R,\epsilon}^p
-(\psi_{R,\epsilon})_t\Big)=\\
&-(\frac2{p-1})\frac {d-1}{2r}(r-\epsilon)(R-r)(R+\epsilon-2r)
(1+r)^{\frac2{p-1}}(1+\frac{\epsilon^l}{r^l}R^\frac2{p-1})\\
&+(\frac2{p-1})\frac {d-1}{2r}((r-\epsilon)(R-r))^2
(1+r)^{\frac2{p-1}-1}(1+\frac{\epsilon^l}{r^l}R^\frac2{p-1})\\
&-l\frac {d-1}{2r}((r-\epsilon)(R-r))^2
(1+r)^{\frac2{p-1}}\frac{\epsilon^l} {r^{l+1}}R^\frac2{p-1}\\
&+(\beta-\gamma)((r-\epsilon)(R-r))^2(1+r)^{\frac2{p-1}}
(1+\frac{\epsilon^l}{r^l}R^\frac2{p-1})\\
&-(1+r)^{\frac{2p}{p-1}}(1+\frac{\epsilon^l}{r^l}R^\frac2{p-1})^p
\exp((p-1)\gamma(t+1)).
\end{aligned}
\end{equation}
We will   show that for all sufficiently large $R$ and sufficiently small
$\epsilon$, and
for sufficiently large $\gamma$ and sufficiently small $l$,
independent of those $R$ and
$\epsilon$, the sum of the right hand sides of \eqref{second}
and \eqref{other} is non-positive.
This will then prove \eqref{supersolution}.

We will denote the seven terms on the right hand side of
\eqref{second} by $J_1-J_7$, and the five terms on the right hand
side of \eqref{other} by $I_1-I_5$.
Note that the terms
that are positive are $J_1, J_2, J_4, J_5, J_7$ and $I_2$.
(Since $\beta$ is bounded from above, $I_4$ is negative for $\gamma$
sufficiently large.)
In what follows, $M$ will denote a positive number that can be made
as large as one desires  by choosing $\gamma$ sufficiently large.
Consider first those $r$ satisfying $r\ge cR$, where $c$ is a fixed positive
number. For $r$ in this range, we have
$|I_5|\ge MR^{\frac2{p-1}+2}(1+\epsilon^l R^{\frac2{p-1}-l})$,
It is easy to see that for $M$ sufficiently large,
 $|I_5|$ dominates each of the positive terms,
uniformly over large $R$ and small $\epsilon$,
and thus (since $M$ can be made arbitrarily large) also the sum of
all of the positive terms. Now consider those $r$ for which
$\delta_0\le r\le C$, for some constants $0<\delta_0<C$.
For $r$ in this range, we have
$|I_4|\ge MR^2(1+\epsilon^l R^{\frac2{p-1}})$, and it is easy to see
that for $M$ sufficiently large, $|I_4|$ dominates
each of the positive terms, uniformly over large $R$ and small $\epsilon$,
and thus, also the sum of all of the positive terms.
One can also show that the transition from $r$ of order unity to
$r$ of order $R$ causes no problem. Thus, we conclude that for any fixed
$\delta_0>0$, for all $l\in(0,1]$ and $\gamma$
sufficiently large, the sum of the  right hand sides of
\eqref{second} and \eqref{other} is negative for all large $R$
and small $\epsilon$.

We now turn to the case
that  $\epsilon\le r\le \delta_0$.
(Note that at $r=\epsilon$, all the terms vanish except
$J_1$ and $I_5$.  Using the fact that
$\frac2{p-1}+2=\frac{2p}{p-1}$,
it is easy to see that for sufficiently large
$\gamma$, $|I_5(\epsilon)|$ dominates $J_1(\epsilon)$, uniformly over all
large $R$ and small $\epsilon$. However, when
$r$ is small, but on an  order larger  than $\epsilon$, the analysis becomes
a lot more involved.)
In the sequel, whenever we say that a condition holds for $\gamma$
or $M$  sufficiently large, or for $l$ sufficiently small,
we mean  that it holds independent of $R$ and $\epsilon$.

Clearly, $J_5\le |I_4|$
if $\gamma$ is sufficiently large.
We now show that for $\gamma$ sufficiently
large,  $J_2\le |I_4|+|I_5|$, for $\epsilon\le r\le\delta_0$.
(We are reusing $|I_4|$ here. Later we will reuse $|I_5|$. This is
permissible because $\gamma$
can be chosen as large as we like.)
To show this inequality, it suffices to show that for $M$ sufficiently large,
\begin{equation}\label{J_2}
(r-\epsilon)R
\le M(r-\epsilon)^2R^2
+M(1+\frac{\epsilon^l}{r^{l+1}}R^\frac2{p-1})^{p-1},\
\text{for}\ r\in[\epsilon,\delta_0].
\end{equation}
A trivial calculation shows that the left hand side of \eqref{J_2}
is less than the first term on the right hand side if
$r\ge\epsilon+\frac1{RM}$. If $r\in[\epsilon,\epsilon+\frac1{RM}]$,
then the left hand side of \eqref{J_2} is less than
or equal to $\frac1M$ while the second term on the right hand side
is greater than $M$. We conclude that \eqref{J_2}
holds with $M\ge1$.

Since $I_2$ has the factor $(r-\epsilon)^2$,
while $I_1$ has the factor $(r-\epsilon)$,
and since $\frac{R-r}{R+\epsilon-2r}$ can be made arbitrarily close
to one by choosing $R$ sufficiently large,
it follows that  for any $\eta>0$, we can guarantee that
$I_2\le\eta|I_1|$, for $r\in[\epsilon,\delta_0]$,
if we choose $\delta_0$ sufficiently small.
Note that $J_4\le\frac {2l}{d-1}|I_1|$. Thus,
given any $\zeta>0$, if
we choose $\delta_0$
and $l$ sufficiently small, we will  have
$I_2+J_4\le\zeta |I_1|$.

To complete the proof, we will show that
\begin{equation}\label{J_1J_7I_4I_5}
J_1+J_7+(1-\zeta)I_1+I_4+I_5\le 0, \ \text{for}\ r\in[\epsilon,\delta_0],
\end{equation}
for sufficiently large $\gamma$ and sufficiently small $l$ and  $\delta_0$,
uniformly over large $R$ and small $\epsilon$.
By the assumption on $\beta$, there exists an $\eta_0>0$ such
that if $\delta_0$ is chosen sufficiently small, then
\begin{equation}\label{I_4}
\begin{aligned}
&I_4\le(\beta_0-\eta_0)(1-\frac\epsilon r)^2(R-r)^2
(1+r)^{\frac2{p-1}}(1+\frac{\epsilon^l}{r^l}R^\frac2{p-1}),\\
&\text{for}\ r\in[\epsilon,\delta_0],
\end{aligned}
\end{equation}
where $\beta_0$ is as in the statement of the theorem.
Since $\frac{R+\epsilon-2r}{R-r}$ can be made arbitrarily close to one
by choosing  $R$ sufficiently large, we also have
\begin{equation}\label{I_1}
I_1\le-(\frac{d-1}{p-1}-\frac{\eta_0}2)(1-\frac\epsilon r) (R-r)^2
(1+r)^{\frac2{p-1}}(1+\frac{\epsilon^l}{r^l}R^\frac2{p-1}),
\ \text{for}\ r\in[\epsilon,\delta_0].
\end{equation}
Since
\[
\begin{aligned}
&J_7=\frac{l(l+1)}2((r-\epsilon)(R-r))^2(1+r)^\frac2{p-1}
\frac{\epsilon^l}{r^{l+2}}R^\frac2{p-1}\\
&\le
\frac{l(l+1)}2(R-r)^2(1+r)^\frac2{p-1}
\frac{\epsilon^l}{r^l}R^\frac2{p-1},
\end{aligned}
\]
and since $\frac{R-r}{R+\epsilon-2r}$ can be made arbitrarily close to  one
by choosing $R$ sufficiently large,
it follows that for any $\tau>0$, we can choose $l$ sufficiently small
so that
\begin{equation}\label{J_7J_1}
J_7\le\tau J_1, \ r\in[\epsilon,\delta_0],
\end{equation}
uniformly over large $R$ and small $\epsilon$.
Using \eqref{I_4}-\eqref{J_7J_1} along with the fact that
\[
J_1\le\frac{p+1}{(p-1)^2}(R-r)^2(1+r)^\frac2{p-1}
(1+\frac{\epsilon^l}{r^l}R^\frac2{p-1})
\]
and that $\beta_0=\frac{d(p-1)-2p}{(p-1)^2}$,
it follows that
\begin{equation}\label{J_1J_7I_4}
\begin{aligned}
&J_1+J_7+(1-\zeta)I_1+I_4\le
\Big(\frac{\tau(p+1)+\zeta(d-1)(p-1)}{(p-1)^2}-\eta_0(\frac12+\frac\zeta2)
+C\frac\epsilon r\Big)\times\\
&(R-r)^2(1+r)^\frac2{p-1}
(1+\frac{\epsilon^l}{r^l}R^\frac2{p-1}), \ \text{for}\ r\in[\epsilon,\delta_0],
\end{aligned}
\end{equation}
for some $C>0$, uniformly over large $R$ and small $\epsilon$.
By picking $\tau$ and $\zeta$ sufficiently small,
we have $\frac{\tau(p+1)+\zeta(d-1)(p-1)}{(p-1)^2}-\eta_0(\frac12+\frac\zeta2)<0$.
Thus, in order to show \eqref{J_1J_7I_4I_5},
it suffices to show that
\begin{equation}\label{final}
\frac\epsilon r R^2\le M(1+\frac{\epsilon^l}{r^l}R^\frac2{p-1})^{p-1},
\ r\in[\epsilon,\delta_0],
\end{equation}
for sufficiently large $M$.
But since $l(p-1)\le1$, the right hand side of \eqref{final}
is greater or equal to $M\frac\epsilon r R^2$.

\end{proof}

\noindent \bf Acknowledgement.\rm\ R. P. thanks  Robert Adler for a helpful
conversation.

\end{document}